\documentclass[12pt]{article}
\usepackage{amssymb,amsmath,enumerate,latexsym, graphicx, epsfig}
\usepackage[latin1]{inputenc}

\def \C {{\Bbb C\,}}
\def\R {{\Bbb R}}

\def \Z {{\Bbb Z}}

\def \esf {{\Bbb S}}

\def\l{{\lambda}} 
\def\a{{\alpha}} 
\def\be{{\beta}}
\def\t{{\theta}} 
\def\g{{\gamma}}

\newtheorem{theorem}{Theorem}
\newtheorem{proposition}{Proposition}
\newtheorem{definition}{Definition}
\newtheorem{lemma}{Lemma}

\newtheorem{remark}{Remark}

\newcommand{\myskip}[1]{}

\begin{document}
\title{A Jenkins-Serrin problem on the strip}
\author{M. Magdalena Rodr\'\i guez\thanks{Research partially supported by grants from R\'egion
  Ile-de-France and a MEC/FEDER grant no.  MTM2004-02746.}}
\date{ }

\maketitle

\noindent {\sc Abstract.} {\footnotesize We describe the family of minimal graphs on strips with boundary
  values $\pm\infty$ disposed alternately on edges of length one,
  and whose conjugate graphs are contained in horizontal slabs of width
  one in $\R^3$.
  We can obtain as limits of such graphs the helicoid, all the doubly
  periodic Scherk minimal surfaces and the singly periodic Scherk
  minimal surface of angle $\pi/2$.}

\section{Introduction}
\label{secintrod}
Karcher~\cite{ka4,ka6} constructed a class of doubly periodic minimal
surfaces, called {\it toroidal halfplane layers}, from minimal graphs,
by extending such graphs by symmetries. More precisely, he considered
the solution to the minimal graph equation on a rectangle with
boundary values~$0$ on the longer edges and $+\infty$ on the shorter
ones; and he extended such a minimal graph to a whole strip by
rotating it an angle $\pi$ about the straight segments corresponding
to the boundary values $0$ (see the upper picture on
Fig.~\ref{grafos}).  The toroidal halfplane layer is obtained from this
Jenkins-Serrin graph on the strip by considering the $\pi$-rotation
about the vertical straight lines on its boundary. Such a doubly
periodic example is denoted by $M_{\t,\frac{\pi}{2},\frac{\pi}{2}}$
in~\cite{mrod1}.  Indeed, this is a particular case in the
3-parametric family of KMR examples $M_{\t,\a,\be}$, with
$\t\in(0,\frac{\pi}{2})$, $\a\in(\frac{- \pi}{2},\frac{\pi}{2}]$,
$\be\in[0,\pi)$ and $(\a,\be)\neq(0,\t)$, examples which have been
classified in~\cite{PeRoTra1} as the only properly embedded, doubly
periodic minimal surfaces with parallel ends and genus one in the
quotient.  Similarly to the construction of
$M_{\t,\frac{\pi}{2},\frac{\pi}{2}}$, Karcher obtained the KMR example
$M_{\t,0,\frac{\pi}{2}}$ by considering the solution to the
Jenkins-Serrin problem on a rectangle with boundary values $0$ on its
longer edges and $+\infty,-\infty$ on its shorter ones (see
Fig.~\ref{grafos}, down).  He also described a continuous
deformation from $M_{\t,\frac{\pi}{2},\frac{\pi}{2}}$ to
$M_{\t,0,\frac{\pi}{2}}$, which corresponds to the surfaces denoted by
$M_{\t,\a,\frac{\pi}{2}}$ in~\cite{mrod1}, with $\a\in[0,\frac{\pi}{2}]$,
and pointed out that the intermediate surfaces did not have enough
symmetries to construct them as Jenkins-Serrin graphs.

We prove that it is possible to construct each
$M_{\t,\a,\frac{\pi}{2}}$, with $\a\in(\frac{- \pi}{2},\frac{\pi}{2}]$,
from a Jenkins-Serrin graph on a parallelogram $\mathcal{P}$ with
boundary values $+\infty$ on its shorter edges and bounded data
$f_1,f_2$ on its longer ones, and this graph can be extended to a
Jenkins-Serrin graph on the strip (see the middle picture in
Fig.~\ref{grafos}). In this case, such an extension does not consist
of a rotation about a straight line, but of the composition of the
reflection symmetry across the plane containing the parallelogram
$\mathcal{P}$ and the translation by the shorter edges on
$\partial\mathcal{P}$. In particular, it must hold $f_1=-f_2$.
Recently, Mazet~\cite{mazet4} has recently constructed, in a
theoretical way, these Jenkins-Serrin graphs on the strip.

Given $h>0$ and $a\in(\frac{-1}{2},\frac{1}{2}]$, consider
$p_n=(n-a,0,-h)$ and $q_n=(n+a,0,h)$, for every $n\in\Z$. We define
the strip $S(h,a)=\{(x_1,0,x_3)\ |\ -h<x_3<h\}$ and mark its boundary
straight lines by $+\infty$ on the straight segments $(p_{2k},
p_{2k+1})$, $(q_{2k}, q_{2k+1})$ and $-\infty$ on $(p_{2k-1}, p_{2k})$,
$(q_{2k-1}, q_{2k})$.
Note that we do not consider $S(h,\frac{-1}{2})$ because it coincides
with $S(h,\frac{1}{2})$.

\begin{definition}\label{def}
  We will say that a minimal graph defined on $S(h,a)$ {\it solves the
    Jenkins-Serrin problem on $S(h,a)$} if its boundary values are
  $\pm\infty$ as prescribed above on each unitary segment
  $(p_n,p_{n+1}),(q_n,q_{n+1})\subset\partial S(h,a)$.
\end{definition}

We know from~\cite{jes1} that, in order to solve the Jenkins-Serrin
problem on $S(h,a)$, it must be satisfied $|q_0-p_0| > 1$; this is,
$a^2+h^2>\frac{1}{4}$.  We define the collection of marked strips
\[
\textstyle{\mathcal{S}=\left\{S(h,a)\ |\ h>0\ \mbox{and }
  a\in\left(\frac{-1}{2},\frac{1}{2}\right]\ \mbox{satisfy }
  a^2+h^2>\frac{1}{4}
\right\}.}
\]

\begin{theorem}
  \label{mainth}
  For every marked strip $S(h,a)\in\mathcal{S}$, there exist
  $\t\in(0,\frac{\pi}{2})$ and $\a\in(\frac{- \pi}{2},\frac{\pi}{2}]$
  such that a piece of the KMR example $M_{\t,\a,\frac{\pi}{2}}$
  solves the Jenkins-Serrin problem on $S(h,a)$.  Moreover, if a
  minimal graph $M$ solves the Jenkins-Serrin problem on some
  $S(h,a)\in\mathcal{S}$ and its conjugate surface is contained in the
  slab $\{(x_1,x_2,x_3)\ |\ 0<x_2<1\}$ up to a translation, then $M$
  must be a piece of a KMR example $M_{\t,\a,\frac{\pi}{2}}$.
\end{theorem}

The author would like to thank Martin Traizet, who motivated
this work during her stay in the University of Fran\c{c}ois Rabelais at Tours.

\section{The KMR examples $M_{\t,\a,\frac{\pi}{2}}$}
\label{secejemplos}

We know~\cite{PeRoTra1} that the space of doubly periodic minimal
surfaces in $\R^3$ with parallel ends and genus one in the quotient
coincides with the family of KMR examples
$\left\{M_{\t,\a,\be}\ |\ \t\in(0,\frac{\pi}{2}),
\a\in(\frac{- \pi}{2},\frac{\pi}{2}], \be\in[0,\pi),
(\a,\be)\neq(0,\t)\right\}$, which has been studied in detail and
classified in~\cite{mrod1} (we will keep the notation introduced
there).
We do not consider the example $M_{\t,\frac{- \pi}{2},\be}$
because it coincides with $M_{\t,\frac{\pi}{2},\be}$,
for every $\t,\be$.
Here we sketch some properties of the subfamily
$\{M_{\t,\a,\frac{\pi}{2}}\}_{\t,\a}$.

Given $\t \in(0,\frac{\pi}{2})$ and $\a \in[0,\frac{\pi}{2}]$,
the minimal surface $M_{\t,\a,\frac{\pi}{2}}$ is determined by the Weierstrass data
\[
\textstyle{g(z,w)=-i+\frac{2}{e^{i\a}z-i}\qquad\mbox{and }\qquad dh=\mu\frac{dz}{w},\quad
\mu\in\R-\{0\}}
\]
(here $g$ is the Gauss map of $M_{\t,\a,\frac{\pi}{2}}$ and $dh$ is
its height differential),
defined on the rectangular torus $\Sigma_\t= \left\{(z,w)\in
  \overline{\C}^2\ |\ w^2=(z^2+\l_\t^2)(z^2+\l_\t^{-2})\right\}$, where
$\l_\t=\cot\frac{\t}{2}$.  The ends of $M_{\t,\a,\frac{\pi}{2}}$, which
are horizontal and of Scherk-type, correspond to the zeroes $A',A'''$
and poles $A,A''$ of~$g$ (i.e. those points with $z=-i e^{-i\a}$ and
$z=i e^{-i\a}$, respectively).  And the Gauss map $g$ of
$M_{\t,\a,\frac{\pi}{2}}$ has four branch points on $\Sigma_\t$:
$D=(-i \l_\t,0)$, $D'=(i \l_\t,0)$,
$D''=(\frac{i}{\l_\t},0)$ and $D'''=(\frac{-i}{\l_\t},0)$.

\begin{figure}\begin{center}
\epsfysize=4cm \epsffile{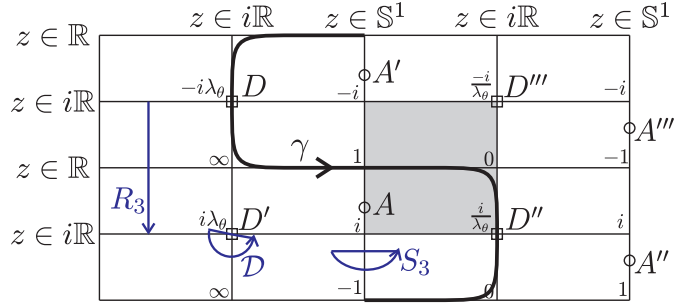}
\end{center}
  \caption{The torus $\widetilde\Sigma_\t$. The value appearing at
    each intersection point between a horizontal and a vertical line
    refers to the value of the $z$-map at the corresponding point.}
\label{torus}
\end{figure}

The multivalued, doubly periodic map $z:\Sigma_\t\to\overline\C$ is
used in~\cite{mrod1} to describe a conformal model of $\Sigma_\t$ as a
quotient of the plane by two orthogonal translations $l_1,l_2$. One of
the advantages is that we can read directly the $z$-values in this
model.  A fundamental domain in $\C$ of the action of the group
generated by $l_1,l_2$ is the parallelogram $\widetilde\Sigma_\t$
represented in Fig.~\ref{torus}.  Each vertical line on
$\widetilde\Sigma_\t$ corresponds to a horizontal level section of
$M_{\t,\a,\frac{\pi}{2}}$ (i.e. a set $x_3^{-1}(constant)$, where
$x_3=\Re\int dh$ on $M_{\t,\a,\frac{\pi}{2}}$). The curve $\g$ drawn
in Fig.~\ref{torus} represents a homology class in
$\Sigma_\t-\{A,A',A'',A'''\}$ with vanishing period.  Since the
periods of $M_{\t,\a,\frac{\pi}{2}}$ at its ends are
\[
\textstyle{\mbox{Per}_A=\mbox{Per}_{A'}=-\mbox{Per}_{A''}=-\mbox{Per}_{A'''}=
  \Big(\mu\,\frac{\pi\,\sin\t}{\sqrt{1-\sin^2\t\,\cos^2\a}},0,0\Big)}\,
,
\]
we conclude that every vertical line in $\widetilde\Sigma_\t$
corresponds to a curve in $\Sigma_\t$ with period $\pm\mbox{Per}_A$. We
fix $\mu$ so that $P=\mbox{Per}_A=(2,0,0)$.

The flux vectors of $M_{\t,\a,\frac{\pi}{2}}$ at its ends are
$\mbox{Fl}_A=-\mbox{Fl}_{A'}=-\mbox{Fl}_{A''}=\mbox{Fl}_{A'''}=
\Big(0,-2,0\Big)$.  Thus we say that $A,A'''$ (resp. $A',A''$) are
{\it left ends} (resp. {\it right ends}).

If we denote by $\widetilde\g\subset\Sigma_\t$ the curve which
corresponds in $\widetilde\Sigma_\t$ to the horizontal line passing
through $D,D'''$, then the flux of $M_{\t,\a,\frac{\pi}{2}}$ along
$\widetilde\g$ equals $-\mbox{Fl}_A$, and the period of
$M_{\t,\a,\frac{\pi}{2}}$ along $\widetilde\g$ can be written as
$T=(T_1,0,T_3)$, with $T_3\neq 0$.  In particular, $T$ is never
horizontal, and $M_{\t,\a,\frac{\pi}{2}}$ is a doubly periodic minimal
surface with period lattice generated by $P,T$.

For every $\t \in(0,\frac{\pi}{2})$ and $\a \in[0,\frac{\pi}{2}]$, we
can similarly define the surface $M_{\t,-\a,\frac{\pi}{2}}$ which
coincides with the reflected image of $M_{\t,\a,\frac{\pi}{2}}$ with
respect to a plane orthogonal to the $x_1$-axis.
Finally, recall from~\cite{mrod1} that the conjugate surface of
$M_{\t,\a,\frac{\pi}{2}}$ coincides (up to normalization) with the
KMR example $M_{\frac{\pi}{2}-\t,\a,0}$, and its periods
(resp. flux vectors) at the ends point to the $x_2$-direction (resp.
$x_1$-direction).

\subsection{Isometries of $M_{\t,\a,\frac{\pi}{2}}$}

The surface $M_{\t,\a,\frac{\pi}{2}}$ has four horizontal straight
lines traveling from left to right ends.  The $\pi$-rotation about any
of those straight lines induce the same isometry $S_3$ of
$M_{\t,\a,\frac{\pi}{2}}$, which corresponds to a symmetry of
$\widetilde\Sigma_\t$ across any of the two vertical lines passing
through the ends.

Another isometry of $M_{\t,\a,\frac{\pi}{2}}$, denoted by
$\mathcal{D}$, is induced by the deck transformation, and corresponds
to the central
symmetry across any of the four branch points of $g$ in
either $\R^3$ or $\widetilde\Sigma_\t$.

The isometry group of $M_{\t,\a,\frac{\pi}{2}}$, which is isomorphic
to $(\Z/2\Z)^3$, is generated by $S_3,\mathcal{D}$ and $R_3$, where
$R_3$ corresponds to the composition of a reflection symmetry across
the plane orthogonal to the $x_2$-axis containing the four branch
points of $g$, with a translation by $(1,0,0)$. The isometry $R_3$
corresponds in $\widetilde\Sigma_\t$ to the translation by half a
vertical period, see Fig.~\ref{torus}.

When $\a=0,\pi/2$, the isometry group of $M_{\t,\a,\frac{\pi}{2}}$ is
richer (it is isomorphic to $(\Z/2\Z)^4$), but we
will not use this fact along this work.
This is the lack of isometries that Karcher referred to for the intermediate
surfaces $M_{\t,\a,\frac{\pi}{2}}$, $0<\a<\frac{\pi}{2}$.

\section{$M_{\t,\a,\frac{\pi}{2}}$ as a graph over the $x_1
  x_3$-plane: $\widetilde S_{\t,\a,\frac{\pi}{2}}$}
Consider the rectangular domain in $\widetilde\Sigma_\t$ on the right of the
middle vertical line.  It corresponds to a piece of
$M_{\t,\a,\frac{\pi}{2}}$ (in fact, we know by $S_3$ that it is a half
of $M_{\t,\a,\frac{\pi}{2}}$), which is a noncompact, singly periodic
minimal annulus bounded by four horizontal straight lines.  We consider
a component $\widetilde S_{\t,\a,\frac{\pi}{2}}$ of the lifting of
this annulus to $\R^3$, and call $S_{\t,\a,\frac{\pi}{2}}$ a
fundamental domain of $\widetilde S_{\t,\a,\frac{\pi}{2}}$ (see
Fig.~\ref{grafos}).

We can assume that $D''$ lies at the origin of $\R^3$ and
$R_3,\mathcal{D}$ are respectively given by the restrictions to
$M_{\t,\a,\frac{\pi}{2}}$ of
\[
(x_1,x_2,x_3)\mapsto(x_1+1,-x_2,x_3) ,\qquad
(x_1,x_2,x_3)\mapsto(-x_1,-x_2,-x_3) .
\]
Take $h>0$ so that the four horizontal straight lines on the boundary
of $S_{\t,\a,\frac{\pi}{2}}$ lie in $\{x_3=\pm h\}$.  Hence both
$S_{\t,\a,\frac{\pi}{2}}$ and $\widetilde S_{\t,\a,\frac{\pi}{2}}$ are
contained in the horizontal slab $\{(x_1,x_2,x_3)\ |\ -h<x_3<h\}$.
Moreover, the horizontal level sections of $S_{\t,\a,\frac{\pi}{2}}$
(which correspond to the vertical lines of $\widetilde\Sigma_\t$ on
the right of the middle vertical line) have period~$P=(2,0,0)$, up to
sign.  Hence $\widetilde S_{\t,\a,\frac{\pi}{2}}$ projects
orthogonally in the $x_2$-direction onto the whole strip
$\mathcal{B}={\{(x_1,0,x_3)\ |\ -h<x_3<h\}}$. Finally, let $\Pi:\widetilde
  S_{\t,\a,\frac{\pi}{2}}\to \mathcal{B}$ be the orthogonal projection
in the $x_2$-direction, $\Pi(p)=(x_1(p),0,x_3(p))$.

\begin{proposition}
  \label{propos}
  The surface $\widetilde S_{\t,\a,\frac{\pi}{2}}$ solves the Jenkins-Serrin
  problem on $S(h,a)$, for some $a\in(\frac{-1}{2},\frac{1}{2}]$.
\end{proposition}
\begin{proof}
  Firstly assume $\widetilde S_{\t,\a,\frac{\pi}{2}}$ is a graph over
  the strip $\mathcal{B}$, $u:\mathcal{B}\to\widetilde
  S_{\t,\a,\frac{\pi}{2}}$.  Recall that $M_{\t,\a,\frac{\pi}{2}}$ has
  horizontal Scherk-type ends with period $(2,0,0)$ and that we obtain
  a fundamental domain of $M_{\t,\a,\frac{\pi}{2}}$ by rotating
  $S_{\t,\a,\frac{\pi}{2}}$ about one of the four straight lines in
  $\partial S_{\t,\a,\frac{\pi}{2}}$. Hence the boundary of
  $\widetilde S_{\t,\a,\frac{\pi}{2}}$ consists of straight lines
  whose orthogonal projection in the $x_2$-direction is formed by two
  rows of equally spaced points, which we can denote by
  $p_n=(n-a,0,-h)$, $q_n=(n+a,0,h)$, for $n\in\Z$ and some
  $a\in(\frac{-1}{2},\frac{1}{2}]$, in such a way that $u$ diverges to
  $+\infty$ when we approach $(p_{2k},p_{2k+1}),(q_{2k}, q_{2k+1})$
  and diverges to~$-\infty$ when we approach $(p_{2k-1},
  p_{2k}),(q_{2k-1}, q_{2k})$ within $\mathcal{B}$, for every
  $k\in\Z$. This is, $\widetilde S_{\t,\a,\frac{\pi}{2}}$ solves the
  Jenkins-Serrin problem on $S(h,a)$, see Definition~\ref{def}.
  Therefore, to conclude Proposition~\ref{propos} it suffices to prove
  that $\widetilde S_{\t,\a,\frac{\pi}{2}}$ is a graph
  over~$\mathcal{B}$.

\begin{figure}
\begin{center}
\epsfysize=9cm \epsffile{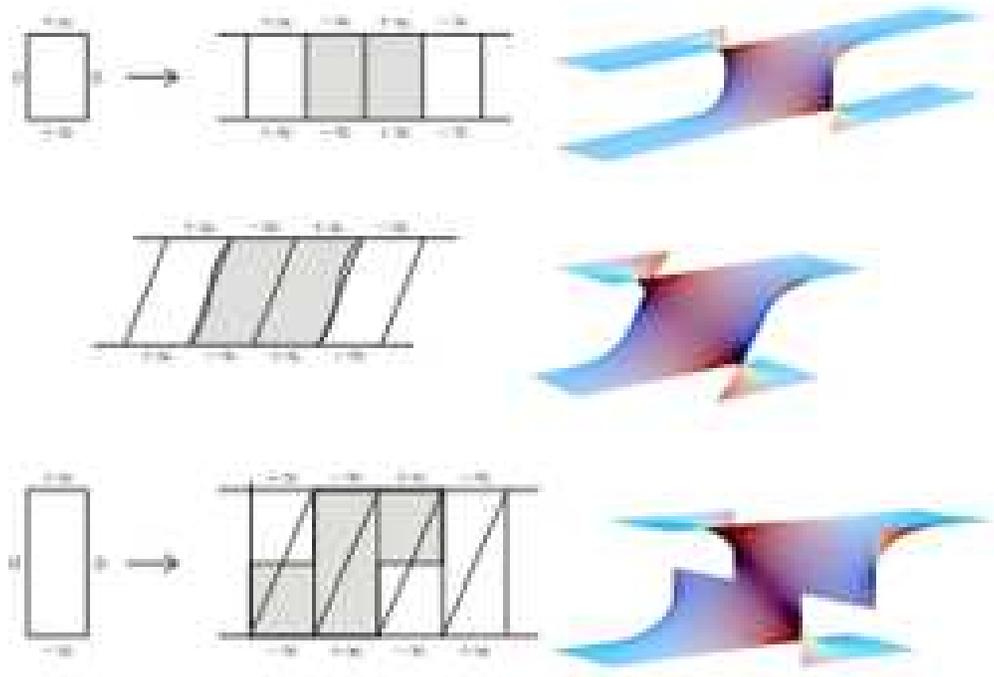}
\end{center}
\caption{Construction of the graphs
  $S_{\frac{\pi}{4},0,\frac{\pi}{2}}$ (top) and
  $S_{\frac{\pi}{4},\frac{\pi}{2},\frac{\pi}{2}}$ (bottom). And the
  intermediate graph $S_{\frac{\pi}{4},\frac{\pi}{4},\frac{\pi}{2}}$
  (center).}
\label{grafos}
\end{figure}

  Denote by $\mathcal{R}$ the piece of $S_{\t,\a,\frac{\pi}{2}}$ which
  corresponds to the region of $\widetilde\Sigma_\t$ shadowed in
  Fig.~\ref{torus}; this is, the rectangle of $\widetilde\Sigma_\t$ on
  the right (resp. left) of the vertical line passing through $A,A'$
  (resp.  $D'',D'''$) and above (resp. below) the horizontal line
  passing through $D',D''$ (resp. $D,D'''$).  The boundary of
  $\mathcal{R}$ consists of a horizontal curve $c_1$ in $\R^3$ joining
  the branch points $D'',D'''$, two curves $c_2,c_3$ from $D'',D'''$,
  respectively, to the horizontal plane $\{(x_1,x_2,x_3)\ |\ x_3=-h\}$
  and either two half straight lines (when $\a\neq 0$) or a straight
  line (when $\a=0$) in $\{x_3=-h\}$.  Since $R_3(c_2)=c_3$, then
  $\Pi(c_3)=\Pi(c_2)+(1,0,0)$.

  Assume $\mathcal{R}$ is a graph over $\mathcal{B}$, and let us prove
  that the same holds for $\widetilde S_{\t,\a,\frac{\pi}{2}}$.
  Suppose by contradiction there exist two points $p,q\in\widetilde
  S_{\t,\a,\frac{\pi}{2}}$ with $\Pi(p)=\Pi(q)$.  In particular,
  $x_1(p)=x_1(q)$. Since $\mathcal{D}(x_1,x_2,x_3)=(-x_1,-x_2,-x_3)$
  and $\mathcal{R}$ is a graph over $\mathcal{B}$, we can assume
  $p\in\mathcal{R}$ and $q\in R_3(\mathcal{R})$.  Let us call $p'$ the
  point in $c_2$ at the same height than $p,q$ (in particular,
  $p',p,q$ correspond to three points in the same vertical line of
  $\widetilde\Sigma_\t$). Hence, by using the isometry $R_3$ and the
  fact that $\mathcal{R}$ is a graph over $\mathcal{B}$, we deduce
  $x_1(p')<x_1(p)<x_1(p')+1<x_1(q)$, a contradiction.

  Therefore, let us prove that $\mathcal{R}$ is a graph over $\mathcal{B}$.
  The spherical image of $\mathcal{R}$ by its Gauss map is contained
  in a quarter of sphere in $\esf^2\cap\{x_2>0\}$, so either
  $\mathcal{R}$ is a graph or a multigraph over $\mathcal{B}$.  The
  following Lemma~\ref{lemma} allows us to conclude that $\mathcal{R}$
  cannot be a multigraph, which finishes Proposition~\ref{propos}.\\ \mbox{}
\end{proof}

\begin{lemma}\label{lemma}
  The restriction of $\Pi$ to $\sigma=c_3\cup c_1\cup c_2$ is one to one.
\end{lemma}
\begin{proof}
  We identify $\sigma$ with its
  corresponding curve in $\Sigma_\t$.  Without loss of generality, we
  can assume that $\sigma$ lies in the same branch of the
  $w$-map (i.e. $w$ is univalent along $\sigma$).  Thus
  we can see $z$ as a parameter on $\sigma$, and so
  $\sigma=\{z= it\ |\ -1<t<1\}$.  In particular, we can
  write the first and third coordinates of $M_{\t,\a,\frac{\pi}{2}}$
  along $\sigma$, denoted by $X_1$ and $X_3$
  respectively, as functions of~$t$.  Since the horizontal level
  sections of $M_{\t,\a,\frac{\pi}{2}}$ correspond to vertical
  segments in $\widetilde\Sigma_\t$, it follows that both
  $X_3|_{c_2},X_3|_{c_3}$ are strictly monotone. Furthermore, the
  restriction of $X_1$ to $c_1=\{z=it\ |\ |t|<\l_\t^{-1}\}$ is also
  strictly monotone because
  \[
  \begin{array}{rl}
    X_1(t)= & \frac{1}{2}\Re\int_{-i\l_\t^{-1}}^{it}\left(\frac{1}{g}-g\right)\, dh\\
    = & \mu\int_{-\l_\t^{-1}}^t \frac{1-s^4}{\left(1-2
          s^2\cos(2\a)+s^4\right)\sqrt{(\l_\t^2-s^2)(\l_\t^{-2}-s^2)}}\, ds .
  \end{array}
  \]
  Since the $\Pi$-projections of $c_1,c_2,c_3$ are separately embedded
  and only intersect at the common extrema, we conclude Lemma~\ref{lemma}.
  \mbox{}\
\end{proof}

\begin{remark}\label{rem}
  Recall that the period lattice of $M_{\t,\a,\frac{\pi}{2}}$ is
  generated by $P=(2,0,0)$ and $T=(T_1,0,T_3)$, $T_3\neq 0$.  Then
  $h=\frac{1}{4}|T_3|$ and $a=\frac{1}{4}|T_1|$ in
  Proposition~\ref{propos}. In particular, it must hold
  $T_1^2+T_3^2>4$.
\end{remark}

\section{Limit graphs of $\widetilde S_{\t,\a,\frac{\pi}{2}}$}
\label{seclimits}

\begin{figure}
\begin{center}
\epsfysize=5cm \epsffile{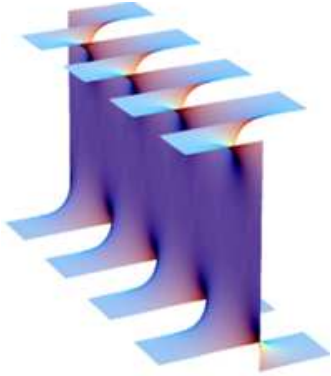}
\end{center}
\caption{The Jenkins-Serrin graph $\widetilde
  S_{\frac{\pi}{200},0,\frac{\pi}{2}}$, close to the singly periodic
  Scherk limit.}
\label{Scherk1p}
\end{figure}

We know~\cite{mrod1} that $M_{\t,\a,\frac{\pi}{2}}$ converges to two
singly periodic Scherk minimal surfaces of angle\footnote{We call {\it
    angle} of a singly or doubly periodic Scherk minimal surface to
  the angle between its nonparallel ends.}~$\frac{\pi}{2}$ when $\t\to
0$.  Let us recall how we can see the singly periodic Scherk minimal
surface of angle~$\frac{\pi}{2}$ as a Jenkins-Serrin graph on the
halfplane.  Consider half a strip $\{0\leq x_1\leq 1,\ x_3\geq 0\}$,
with boundary data~$0$ on the vertical half straight lines and
$+\infty$ on the unit straight segment in between. By rotating about
the boundary half lines, we obtain a Jenkins-Serrin graph
$\widetilde{\cal S}_{1p}$ on the halfplane with boundary values
$\pm\infty$ on $\{x_3=0\}$ disposed alternately on unitary edges,
which is half a singly periodic Scherk minimal surface of angle
$\pi/2$ and period $(2,0,0)$.\\ We have proven that $\widetilde
S_{\t,\a,\frac{\pi}{2}}$ is a graph over the marked strip $S(h,a)$,
where ${h=\frac{1}{4}|T_3|}$ and $a=\frac{1}{4}|T_1|$.  Translate
$\widetilde S_{\t,\a,\frac{\pi}{2}}$ by $(a,0,h)$. Then this translated
$\widetilde S_{\t,\a,\frac{\pi}{2}}$ converges to $\widetilde{\cal
  S}_{1p}$, when $\t\to 0$ (see Fig.~\ref{Scherk1p}).  By using the
isometry $\mathcal{D}$, we obtain that the translated $\widetilde
S_{\t,\a,\frac{\pi}{2}}$ by $(-a,0,-h)$ has a similar behavior.  In
particular, when $\t\to 0$, the width of the strip diverges to
$+\infty$ (i.e. ${|T_3|\to+\infty}$).

\begin{figure}
  \begin{center}
\epsfysize=2cm \epsffile{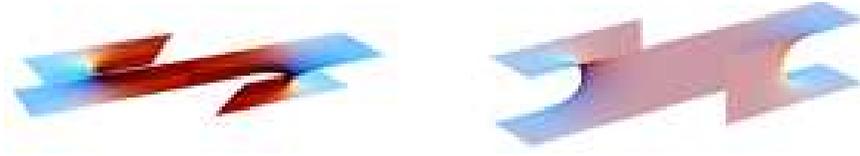}
  \end{center}
  \caption{The Jenkins-Serrin graphs
    $S_{\frac{49\pi}{100},\frac{\pi}{4},\frac{\pi}{2}}$ (left) and
    $S_{\frac{49\pi}{100},\frac{\pi}{2},\frac{\pi}{2}}$ (right), close
    to doubly periodic Scherk minimal surfaces.}
  \label{Scherk2p}
\end{figure}

When $\t\to\frac{\pi}{2}$ and $\a\to\a_\infty\neq 0$,
$M_{\t,\a,\frac{\pi}{2}}$ converges to two doubly periodic Scherk
minimal surfaces of angle $\a_\infty$ and periods of length
one. 
Half a such doubly periodic Scherk example can be seen as a
Jenkins-Serrin graph $\mathcal{S}_{2p}$ on the corresponding rhombus
with alternating boundary data $\pm\infty$.\\
Denote by $\mathcal P_n$ the rhombus of vertices
$p_n,p_{n+1},q_{n+1},q_n$, for every $n\in\Z$, and let $M_n$ be the
piece of $\widetilde S_{\t,\a,\frac{\pi}{2}}$ over $\mathcal P_n$,
translated so that $x_2=0$ in the middle point of $M_n$ (i.e. the
point in $M_n$ which projects onto the middle point of~$\mathcal
P_n$).  For any $k\in\Z$, $M_{2k}$ converges to $\mathcal{S}_{2p}$,
when $\t\to\frac{\pi}{2}$ and $\a\to\a_\infty$ (see
Fig.~\ref{Scherk2p});
and $M_{2k-1}$ converges to the reflected image of $\mathcal{S}_{2p}$
across the $x_1 x_3$-plane.
In this case, $T_1^2+T_3^2\to 4$ and $T_3\not\to 0$.
Moreover, for each $T_{1,\infty},T_{3,\infty}$ with
$T_{1,\infty}^2+T_{3,\infty}^2=4$ and $T_{3,\infty}\neq 0$, there
exists a $\mathcal{S}_{2p}$ which is graph over the parallelogram
determined by $(1,0,0),(T_{1,\infty},0,T_{3,\infty})$; and this
$\mathcal{S}_{2p}$ is obtained as a limit of translated graphs
$S_{\t,\a,\frac{\pi}{2}}$.

\begin{figure}
  \begin{center}
\epsfysize=25mm \epsffile{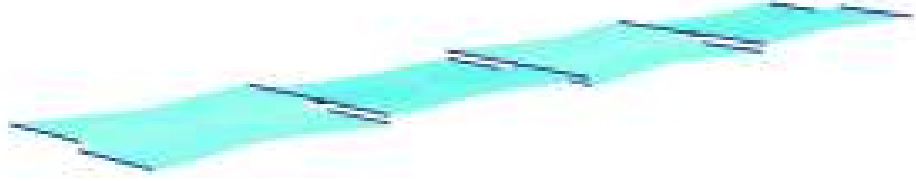}
  \end{center}
  \caption{The Jenkins-Serrin graph $\widetilde
    S_{\frac{49\pi}{100},0,\frac{\pi}{2}}$, close to the helicoid
    limit.}
  \label{Helicoid}
\end{figure}

When $\t\to\frac{\pi}{2}$ but $\a\to 0$, the dilated KMR example
$\frac{1}{\mu}M_{\t,\a,\frac{\pi}{2}}$ converges to two vertical
helicoids spinning oppositely.  Let $\mathcal{H}$ be half a
fundamental domain of the vertical helicoid bounded by two horizontal
straight lines, both projecting vertically onto the same straight line
$\ell\subset\{x_3=0\}$.
Assume $x_1(\ell)=0$ and that the projection of $\partial\mathcal{H}$
in the $x_2$-direction consists of two points at heights $-h$ and $h$.
Thus the interior of $\mathcal{H}$ can be seen as a graph onto the
strip $\{(x_1,0,x_3)\ |\ -h<x_3<h\}$, with boundary data $+\infty$ on
$\{x_1>0,x_2=0,x_3=h\}\cup\{x_1<0,x_2=0,x_3=-h\}$, and $-\infty$ on
$\{x_1>0,x_2=0,x_3=-h\}\cup\{x_1<0,x_2=0,x_3=h\}$.  As
$\t\to\frac{\pi}{2}$ and $\a\to 0$, the suitably translated graphs
$\frac{1}{\mu}S_{\t,\a,\frac{\pi}{2}}$ converge to $\mathcal{H}$ (see
Fig.~\ref{Helicoid}). And different translations of the surfaces
$\frac{1}{\mu}S_{\t,\a,\frac{\pi}{2}}$ converge, when
$\t\to\frac{\pi}{2}$ and $\a\to 0$, to another half a vertical
helicoid spinning oppositely.
In this case, $T_1^2+T_3^2\to 4$ and $T_3\to
0$.

\section{Proof of Theorem~\ref{mainth}}
Denote by $\mathcal{M}$ the family of graphs
\[
\textstyle{
\mathcal{M}=\left\{\widetilde S_{\t,\a,\frac{\pi}{2}}\ \Big|\
\t\in\left(0,\frac{\pi}{2}\right), \a\in\left(\frac{- \pi}{2},\frac{\pi}{2}\right]\right\} .}
\]
Recall that we could have defined the graphs $\widetilde
S_{\t,\frac{- \pi}{2},\frac{\pi}{2}}$ in a similar way, but
$\widetilde S_{\t,\frac{- \pi}{2},\frac{\pi}{2}}=\widetilde S_{\t,\frac{\pi}{2},\frac{\pi}{2}}$.
From the classification of the KMR examples~\cite{mrod1}, we know that
no two surfaces in $\mathcal{M}$ coincide.
This family $\mathcal{M}$ can be naturally endowed with the product
topology given by its parameters $(\t,\a)$.
Furthermore, we know the surfaces obtained by taking limits from graphs in
$\mathcal{M}$ (see Section~\ref{seclimits}).  We deduce that the
boundary $\partial\mathcal{M}$ of $\mathcal{M}$ has two components: an
isolated point $\{\star\}$ corresponding to the singly periodic Scherk
limit $\widetilde S_{1p}$, and a closed curve $\Gamma$ corresponding
to the union of the family of doubly periodic Scherk limits and the
helicoidal limit (recall that the helicoid can be obtained as a limit
surface of doubly periodic Scherk minimal examples).  Hence
$\mathcal{M}$ is topologically a punctured disk $D-\{\star\}$, where
$\Gamma$ is the boundary of the disk~$D$.

Recall the collection of marked strips defined just after
Definition~\ref{def},
\[
\textstyle{
\mathcal{S}=\left\{S(h,a)\ \Big|\ h>0\ \mbox{and }
  a\in\left(\frac{-1}{2},\frac{1}{2}\right]\ \mbox{satisfy }
  a^2+h^2>\frac{1}{4}\right\}.}
\]
Since $S(h,\frac{-1}{2})=S(h,\frac{1}{2})$, the family ${\cal S}$ can
be topologized by the natural map $S(h,a)\in {\cal
  S}\stackrel{H}{\mapsto }(h,a)\in \R^+\times (\R / \Z )$.
Note that the parameter $a$ goes
necessarily to $0$ when $S(h,a)\in\mathcal{S}$ and $h\to+\infty$.
After
identifying $\mathcal{S}$ with its image through $H$, 
we obtain that $\mathcal{S}$ is topologically a punctured disk $D-\{\star\}$, 
and the boundary of
$\mathcal{S}$ consists of two components: the curve $\{(h,a)\ |\
h^2+a^2=\frac{1}{4}\}$, which corresponds to $\Gamma=\partial D$, and
$\{(+\infty,0)\}$, which corresponds to $\{\star\}$.

Proposition~\ref{propos} and Remark~\ref{rem} let us define the continuous map
\[
\begin{array}{rrccc}
  \phi: & \mathcal{M}\equiv D-\{\star\}& \ \longrightarrow\  & \mathcal{S}\equiv
  D-\{\star\} & ,\\
  & \widetilde S_{\t,\a,\frac{\pi}{2}}  & \ \mapsto\ &
  S(\frac{1}{4}|T_3|,\frac{1}{4}|T_1|) &
\end{array}
\]
which can be continuously extended to the boundaries so that $\phi(\star)=\star$
and $\phi(\partial D)=\partial D$, using Section~\ref{seclimits}.

Since the conjugate graph of $\widetilde S_{\t,\a,\frac{\pi}{2}}$ is contained in
$\{(x_1,x_2,x_3)\ |\ 0<x_2<1\}$, then the following lemma implies
that $\phi$ is injective.

\begin{lemma}[Mazet, \cite{mazet3}]
\label{maz}
Let $\Omega$ be a convex polygonal domain with unitary edges, and $M$
be a minimal (vertical) graph on $\Omega$ with boundary data $\pm\infty$ disposed
alternately, and whose conjugate graph lies on a horizontal slab of width one.
Then $M$ is unique up to a vertical translation.
\end{lemma}

It is not difficult to obtain that $\phi$ is onto from the
fact that it is continuous, injective and $\phi(\star)=\star$,
$\phi(\partial D)=\partial D$. This proves the first part in
Theorem~\ref{mainth}.

The uniqueness part in Theorem~\ref{mainth} about the graphs
$\widetilde S_{\t,\a,\frac{\pi}{2}}$ can be also deduced from
Lemma~\ref{maz} as above. This finishes the proof of Theorem~\ref{mainth}.


\end{document}